\documentclass[11pt]{amsart}
\usepackage{amsmath}
\usepackage{amsfonts}
\usepackage{amssymb}
\usepackage{amsthm}
\usepackage{amstext,amscd}

\theoremstyle{plain}
\newtheorem{theorem}{Theorem}[section]

\newtheorem{conjecture}{Conjecture}
\newtheorem{corollary}[theorem]{Corollary}

\newtheorem{definition}[theorem]{Definition}

\newtheorem{lemma}[theorem]{Lemma}

\newtheorem{proposition}[theorem]{Proposition}
\newtheorem*{remark}{Remark}

\numberwithin{equation}{section}
\begin{document}
\title{Volume Conjecture, Regulator and $SL_2(\mathbb{C})$-Character Variety of a Knot}
\author{Weiping Li}
\address{Department of Mathematics\\
         Oklahoma State University\\
         Stillwater, OK 74078}
\email{wli@math.okstate.edu}
\author{Qingxue Wang}
\address{Department of Mathematics\\
         Oklahoma State University\\
         Stillwater, OK 74078}
\email{qinwang@math.okstate.edu}

\begin{abstract}
In this paper, by using the regulator map of Beilinson-Deligne, we
show that the quantization condition posed by Gukov is true for
the $SL_2(\mathbb{C})$ character variety of the hyperbolic knot in
$S^3$. Furthermore, we prove that the corresponding
$\mathbb{C}^{*}$-valued $1$-form is a secondary characteristic
class (Chern-Simons) arising from the vanishing first Chern class
of the flat line bundle over the smooth part of the character
variety, where the flat line bundle is the pullback of the
universal Heisenberg line bundle over $\mathbb{C}^{*}\times
\mathbb{C}^{*}$.

The second part of the paper is to define an algebro-geometric
invariant of $3$-manifolds resulting from the Dehn surgery along a
hyperbolic knot complement in $S^3$. We establish a Casson type
invariant for these $3$-manifolds. In the last section, we
explicitly calculate the character variety of the figure-eight
knot and discuss some applications.
\end{abstract}
\date{April 2, 2006}
\maketitle
\section{Introduction}
It is a very important question in the knot theory to find the
geometric and topological interpretation of the Jones polynomial
of a knot. Observed by Kashaev \cite{Kas}, H. Murakami and J.
Murakami \cite{MM}, the asymptotic rate of (N-colored) Jones
polynomial is related to the volume of the hyperbolic knot
complement. It is known as the Volume conjecture. Following
Witten's $SU(2)$ topological quantum field theory, Gukov
\cite{Guk} proposed a complex version of Chern-Simons theory and
generalized the volume conjecture to a $\mathbb{C}^{*}$-
parametrized version with parameter lying on the zero locus of the
$A$-polynomial of the knot in $S^3$.

In this paper, we prove that the quantization condition posed by
Gukov \cite[Page~597]{Guk} is true for the hyperbolic knots in
$S^3$. The key ingredient of the proof is the construction of the
regulator map of an algebraic curve studied by Beilinson, Bloch,
Deligne and many others. Let $K$ be a hyperbolic knot in $S^3$.
For each irreducible component $Y$ of the zero locus of the
$A$-polynimial $A(l,m)$ of $K$, we show that the symbol
$\{l,m\}\in K_2(\mathbb{C}(Y))$ is a torsion. This element gives
rise to a cohomology class $r(l,m)$ in $H^1(Y_h, \mathbb{C}^{*})$,
where $Y_h$ is some open Riemann surface. The detail is given in
Section $3$. As Deligne noted, $H^1(Y_h, \mathbb{C}^{*})$ is the
group of isomorphism classes of flat line bundles over $Y_h$.
Thus, our class $r(l,m)$ corresponds to a flat line bundle.
Indeed, this line bundle can be constructed explicitly as the
pullback of the universal Heisenberg line bundle over
$\mathbb{C}^{*}\times \mathbb{C}^{*}$, see \cite{Bl,Ram}. We then
derive a $1$-form from this flat line bundle and show that this
\emph{$1$-form} is the Chern-Simons class of the \emph{first}
Chern class $C_1$ of the flat line bundle $r(l,m)$. Note that this
is not the usual Chern-Simons class as a \emph{$3$-form} of the
\emph{second} Chern class for a $3$-dimensional manifold. We also
reformulate the generalized volume conjecture via this $1$-form
Chern-Simons class.

Let $X_0$ be the irreducible component of the character variety of
$K$ which contains the character of the discrete faithful
representation associated to the hyperbolic structure of $K$. By
\cite{CGLS}, $X_0$ is an affine curve. We study in detail the
properties of its image in the character variety of the boundary
(which is a torus) and the corresponding part in the
$A$-polynomial. Using $X_0$, we construct a new $SL_2(\mathbb{C})$
algebro-geometric invariant $\lambda(p,q)$ for the manifolds
obtained by the $(p,q)$ Dehn surgery along the knot complement.
Roughly speaking, our invariant $\lambda(p,q)$ is counting the
intersection multiplicity of $X_0$ with another affine curve in
the character variety of the boundary. It does not require the
non-sufficiently large condition in \cite{Cu1}. Via our invariant
and Culler-Shalen norm, we obtain an upper bound for the number of
ideal points (with multiplicity) which are zeroes of the function
$f_{\gamma}$. The definition of $f_{\gamma}$ is in Section $4$ and
its degree is the Culler-Shalen norm.

Our point of view is that the component $X_0$ should contain a lot
of topological and geometric information about the knot and the
Dehn fillings. Of course, the whole character variety may contain
more information, but at present we still have little information
about other components. It seems that the whole character variety
is complicated in general.

The paper is organized as follows. In section $2$, we introduce
the notations used in the paper. In section $3$, we discuss the
generalized volume conjecture and the regulator of a curve, then
we prove the theorem about the quantization condition. In section
$4$, we define the algebro-geometric invariant $\lambda(p,q)$ of a
hyperbolic knot, and study its properties. In the last section, we
explicitly calculate several invariants of the figure-eight knot,
for instance, its character variety, ideal points, and
Culler-Shalen norm. Some applications of these calculations are
also discussed.
\section{Terminology and Notation}
\subsection{}
Let $K$ be a knot in $S^3$ and $M_K$ its complement. That is,
$M_K=S^3-N_K$ where $N_K$ is the open tubular neighborhood of $K$
in $S^3$. $M_K$ is a compact $3$-manifold with boundary $\partial
M_K=T^2$ a torus. Denote by
$R(M_K)=\text{Hom}(\pi_1(M_K),SL_2(\mathbb{C}))$ and $R(\partial
M_K)=\text{Hom}(\pi_1(\partial M_K),SL_2(\mathbb{C}))$. By
\cite{CS}, they are affine algebraic sets over the complex numbers
$\mathbb{C}$ and so are the corresponding character varieties
$X(M_K)$ and $X(\partial M_K)$. We also have the canonical
surjective morphisms $t:R(M_K)\longrightarrow X(M_K)$ and
$t:R(\partial M_K)\longrightarrow X(\partial M_K)$ which map a
representation to its character. The natural homomorphism $i:
\pi_1(\partial M_K)\longrightarrow \pi_1(M_K)$ induces the
restriction maps $r: X(M_K)\longrightarrow X(\partial M_K)$ and
$r: R(M_K)\longrightarrow R(\partial M_K)$.
\subsection{}
Throughout this paper, for a matrix $A\in SL_2(\mathbb{C})$, we
denote by $\sigma(A)$ its trace.

\subsection{}
Since $\pi_1(\partial M_K)=\mathbb{Z}\oplus \mathbb{Z}$, we shall
fix two oriented simple curves $\mu$ and $\lambda$ as its
generators. They are called the meridian and longitude
respectively. Let $R_D$ be the subvariety of $R(\partial M_K)$
consisting of the diagonal representations. Then $R_D$ is
isomorphic to $\mathbb{C}^{*} \times \mathbb{C}^{*}$. Indeed, for
$\rho \in R_D$, we obtain
\begin{equation*}
\rho(\lambda)=
  \left[
  \begin{matrix}
    l & 0\\
    0 & l^{-1}
  \end{matrix}
  \right] \, \text{and} \, \,
\rho(\mu)=
  \left[
  \begin{matrix}
    m & 0\\
    0 & m^{-1}
  \end{matrix}
  \right],
\end{equation*}
then we assign the pair $(l,m)$ to $\rho$. Clearly this is an
isomorphism. We shall denote by $t_{D}$ the restriction of the
morphism $t:R(\partial M_K)\longrightarrow X(\partial M_K)$ on
$R_D$.

\subsection{}
Next we recall the definition of the $A$-polynomial of $K$ which
was introduced in \cite{CCGLS}. Denote by $X^{\prime}(M_K)$ the
union of the irreducible components $Y^{\prime}$ of $X(M_K)$ such
that the closure $\overline{r(Y^{\prime})}$ in $X(\partial M_K)$
is $1$-dimensional. For each component $Z^{\prime}$ of
$X^{\prime}(M_K)$, denote by $Z$ the curve
$t_{\Delta}^{-1}(\overline{r(Y^{\prime})})\subset R_D$.  We define
$D_K$ to be the union of the curves $Z$ as $Z^{\prime}$ varies
over all components of $X^{\prime}(M_K)$. Via the above
identification of $R_D$ with $\mathbb{C}^{*} \times
\mathbb{C}^{*}$, $D_K$ is a curve in $\mathbb{C}^{*} \times
\mathbb{C}^{*}$. Now by definition the \emph{$A$-polynomial
$A(l,m)$} of $K$ is the defining polynomial of the closure of
$D_K$ in $\mathbb{C}\times\mathbb{C}$.


From now on, we shall assume that $K$ is a hyperbolic knot. Denote
by $\rho_0: \pi_1(M_K)\longrightarrow PSL_2(\mathbb{C})$ the
discrete, faithful representation corresponding to the hyperbolic
structure on $M_K$. Note that $\rho_0$ can be lifted to a
$SL_2(\mathbb{C})$ representation. Moreover, there are exactly
$|H^{1}(M_K; \mathbb{Z}_2)=\mathbb{Z}_2|=2$ such lifts.

\section{$A$-polynomial, regulator and $K_2$ of a curve}
In this section, we briefly recall the generalized volume
conjecture. Using the regulator map of a curve, we show that the
form $r(l,m)=\xi(l,m)+i \eta(l,m)$ over the $1$-dimensional
character variety $Y_h$ has exact imaginary part and the
rationality of the real part. This provides an affirmative answer
to Gukov's quantization over $Y_h$. Moreover,
$dr(l,m)=\frac{dl}{l}\wedge \frac{dm}{m}=2 \pi i C_1(L)=0$
justifies that $r(l,m)$ is the Chern-Simons class from the first
Chern class $C_1$.

\subsection{}
Let $\overline{D_K}$ be the zero locus of the $A$-polynomial
$A(l,m)$ in $\mathbb{C}^2$. Let $y_0\in D_K$ correspond to the
character of the representation of the hyperbolic structure on
$M_K$ and $m(y_0)=1$. For a path $c$ in $\overline{D_K}$ with the
initial point $y_0$ and endpoint $(l,m)$, the following quantities
are defined in \cite[~(5.2)]{Guk}:
\begin{equation}\label{vol}
   Vol(l,m)=Vol(K)+2\int_{c}[-\log{|l|}d(\arg{m})+\log{|m|}d(\arg{l})],
\end{equation}
\begin{equation}\label{cs}
   CS(l,m)=CS(K)-\frac{1}{\pi^2} \int_{c}[\log{|m|}d
\log{|l|}+(\arg{l})d(\arg{m})],
\end{equation}
where $Vol(K)$ and $CS(K)$ are the volume and the Chern-Simons
invariant of the complete hyperbolic metric on $M_K$.

In \cite[~(5.12)]{Guk}, Gukov proposes the \emph{Generalized
Volume Conjecture}: for a fixed number $a$ and
$m=-\exp({i \pi a})$,\\
\begin{equation}\label{conj}
\lim_{N, k\rightarrow \infty; \frac{N}{k}=a}
\frac{\log{J_{N}(K,e^{2\pi i/k})}}{k}=\frac{1}{2 \pi}(Vol(l,m)+
i2\pi^2 CS(l,m)),
\end{equation}
where $J_{N}(K,q)$ is the $N$-colored Jones polynomial of $K$,
$Vol(l,m)$ and $CS(l,m)$ as in (\ref{vol}) and (\ref{cs}), are the
functions on the zero locus of the A-polynomial of the hyperbolic
knot $K$. Note that for $m=1$ and $a=1$, we get the usual
\emph{Volume Conjecture}:
\begin{equation}
\lim_{N\rightarrow \infty} \frac{\log{|J_{N}(K,e^{2\pi
i/N})|}}{N}=\frac{1}{2 \pi}Vol(K).
\end{equation}

The (generalized) volume conjecture links the Jones invariants of
knots with the topological and geometric invariants arising from
the character variety of the knot complement. The volume
conjecture has received a lot of attentions. But other than few
examples being checked, there is no essential mathematical
evidence to support this interesting conjecture. Understanding
those terms in the generalized volume conjecture (\ref{conj})
would be the first step.\\

For $Vol(l,m)$ in (\ref{vol}), it is understood that it measures
the change of volumes of the representations on the path $c$ in
$\overline{D_K}$, the zero locus of $A$-polynomial of the
hyperbolic knot. See \cite[Sect.~4.5]{CCGLS} and
\cite[Sect.~2]{Dun} for more detail.\\

For $CS(l,m)$ in (\ref{vol}), to our knowledge, it has not been
understood mathematically. It was derived from the point of view
of physics, see \cite[Sect.~3]{Guk}. On the other hand, since
$M_K$ has boundary a torus $T$, by \cite{RSW} and \cite{KK}, its
$3$-form Chern-Simons functional is only well-defined as a section
of a circle bundle over the gauge equivalence classes of $T$. By
\cite[Theorem~3.2, 2.7]{KK}, if $\chi_t$, $t\in [0,1]$ is a path
of characters of $SL_2(\mathbb{C})$ representations of $M_K$ and
$z(t)$ is the Chern-Simons invariant of $\chi_t$, then:
\begin{equation}\label{kkf}
  \begin{aligned}
     z(1)z(0)^{-1}&=\text{exp}(2\pi i(\int_{0}^{1}\alpha
                   \frac{d\beta}{dt}-\beta\frac{d\alpha}{dt}))\\
                  &=\text{exp}(\frac{1}{2\pi
                  i}\int_{0}^{1}(\log{m}\;d \log{l}-log{l}\;d
                  \log{m}))
   \end{aligned}
\end{equation}

where $(\alpha(t), \beta(t))$ is a lift of $\chi_t$ to
$\mathbb{C}^2$, under the $(l,m)$ coordinates,
$\displaystyle{\alpha=\frac{1}{2\pi i}\log{m}} $ and
$\displaystyle{\beta=\frac{1}{2\pi i}\log{l}}$ for a fixed branch
of logarithm.\\

Up to scalar, (\ref{kkf}) and (\ref{cs}) are not the same.
Moreover, the Chern-Simons $3$-form in \cite{KK} is the secondary
class from the second Chern class (a closed $4$-form). The term in
(\ref{cs}) defined in \cite[~(5.6)]{Guk} is a $1$-form. It may be
the secondary class of the first Chern class (a closed $2$-form)
of some line bundle over $\overline{D_K}$.\\

In the following subsections, we show that $dCS(l,m)$ is indeed
arising from the first Chern class of a (universal) line bundle
over the Heisenberg group. Furthermore, we relate both $dVol$ and
$dCS$ to the imaginary and real parts of the secondary
Chern-Simons class respectively. We also give a mathematical proof
of Gukov's quantization statement of the Bohr-Sommerfield
condition by using some torsion element of $K_2$ and the regulator
map.

\subsection{The regulator map of $K_2$}
Let $X$ be a smooth projective curve over $\mathbb{C}$ or a
compact Riemann surface. Let $f$, $g$ be two meromorphic functions
on $X$. Denote by $S(f)$ (resp. $S(g)$) the set of zeros and poles
of $f$ (resp. $g$). Notice that $S(f)\cup S(g)$ is a finite set.
Put $X^{\prime}=X\setminus (S(f)\cup S(g))$.

Following Beilinson \cite{Bei}, see also \cite{De}, we define an
element $r(f,g)\in H^{1}(X^{\prime},\mathbb{C}^{*})$,
equivalently, as an element of
$\text{Hom}(\pi_{1}(X^{\prime}),\mathbb{C}^{*})$: for a loop
$\gamma$ in $X^{\prime}$ with a distinguished base point $t_0\in
X^{\prime}$,
\begin{equation}\label{eq2.1}
r(f,g)(\gamma)=\exp{(\frac{1}{2\pi
i}(\int_{\gamma}\log{f}\;\frac{dg}{g}-\log{g(t_0)}\int_{\gamma}\frac{df}{f}))},
\end{equation}

where the integrals are taken over $\gamma$ beginning at $t_0$.

It is well-known that this definition is independent of the
choices of the base point $t_0$ and the branches of $\log{f}$ and
$\log{g}$. From now on, we shall take
$\log{z}:\mathbb{C}^{*}\rightarrow \mathbb{C}$ with $0\leq
\arg{z}< 2 \pi$. Then it is well-defined, but discontinuous on the
positive real line $[0,+\infty)$ and it is holomorphic on the cut
plane $\mathbb{C}\setminus [0,+\infty)$.

In \cite{De}, Deligne noticed that
$H^{1}(X^{\prime},\mathbb{C}^{*})$ is the group of isomorphism
classes of the line bundles over $X^{\prime}$ with flat
connections. Hence $r(f,g)$ corresponds to such a line bundle with
a flat connection.
\begin{proposition}\label{prop1}
(1) The curvature of the line bundle associated to the class
$r(f,g)$ is
$\frac{df}{f}\wedge \frac{dg}{g}$; \\
(2) $r(f_{1}f_{2},g)=r(f_1,g)\otimes r(f_2,g)$,
$r(f,g)=r(g,f)^{-1}$,
and the Steinberg relation $r(f,1-f)=1$ holds if $f\ne 0$, $f\ne 1$;\\
(3) For $x\in S(f)\cup S(g)$, let $\gamma_{x}$ be a small simple
loop around $x$ in $X^{\prime}$. Then $r(f,g)(\gamma_{x})$ is
equal to the tame symbol $T_{x}(f,g)$ of $f$ and $g$ at $x$.
\end{proposition}

\begin{proof}
See \cite{De}. For the explicit construction of the line bundle
$r(f,g)$, see \cite[Section~4]{Ram} and \cite{Bl} where the proof
of this proposition was also given. The key construction is a
universal Heisenberg line bundle with connection on
$\mathbb{C}^{*}\times \mathbb{C}^{*}$. To prove the Steinberg
relation, the ubiquitous dilogarithm shows up.
\end{proof}

Recall the tame symbol
\[
   T_{x}(f,g):=(-1)^{v_{x}(f)\cdot
            v_{x}(g)}\frac{f^{v_{x}(g)}}{g^{v_{x}(f)}}(x),
\]
where $v_{x}(f)$ (resp. $v_{x}(g)$) is the order of zero or pole
of $f$ (resp. $g$) at $x$.\\

Let $\mathbb{C}(X)$ be the field of meromorphic functions on $X$.
Denote by $\mathbb{C}(X)^{*}$ the set of non-zero meromorphic
functions on $X$. By Matsumoto Theorem \cite{Mil},
\[
  K_2(\mathbb{C}(X))=\frac{\mathbb{C}(X)^{*}\otimes
               \mathbb{C}(X)^{*}}{\langle f\otimes (1-f): f\ne 0,1 \rangle},
\]

where the tensor product is taken over $\mathbb{Z}$, and the
denominator means the subgroup generated by those elements. For
$f$ and $g\in \mathbb{C}(X)^{*}$, we denote by $\{f,g\}$ the
corresponding element in $K_2(\mathbb{C}(X))$.

The part $(2)$ of Proposition \ref{prop1} implies that we have a
homomorphism
\begin{equation}\label{eq2.2}
r: K_2(\mathbb{C}(X))\longrightarrow \underset{S\subset
X(\mathbb{C}):\;\text{finite}}{\varinjlim}H^{1}(X\setminus
S,\mathbb{C}^{*})
\end{equation}

defined by $r(\{f,g\})=r(f,g)$.

\subsection{}
Let $Y$ be an irreducible component of $\overline{D_K}$, the zero
locus of the $A$-polynomial $A(l,m)$. Denote by $\widetilde{Y}$ a
smooth projective model of $Y$. Then their fields of rational
functions are isomorphic,
$\mathbb{C}(Y)\cong\mathbb{C}(\widetilde{Y})$. We have the
following.
\begin{proposition}\label{prop2.1}
The element $\{l,m\}\in K_2(\mathbb{C}(Y))$ is a torsion element.
\end{proposition}

\begin{proof}
By \cite[Proposition~2.2, 4.1]{CCGLS}, there is a finite field
extension $F$ of $\mathbb{C}(Y)$ such that $\{l,m\}\in K_2(F)$ is
of order at most $2$. We have a homomorphism $i:
K_2(\mathbb{C}(Y))\rightarrow K_2(F)$ induced by the inclusion of
$\mathbb{C}(Y)$ into $F$. We also have the transfer map
$t:K_2(F)\rightarrow K_2(\mathbb{C}(Y))$. It is well-known that
the composition $t\circ i$:
\[
  K_2(\mathbb{C}(Y))\rightarrow K_2(F)\rightarrow K_2(\mathbb{C}(Y))
\]
is given by the multiplication of $n=[F:\mathbb{C}(Y)]$, the
degree of the finite extension. Hence
$t(i(\{l,m\})=t(\{l,m\})=n\{l,m\}$. This implies that $\{l,m\}\in
K_2(\mathbb{C}(Y))$ is a torsion and its order divides $2n$.
\end{proof}

Suppose the component $Y$ contains $y_0\in D_K$ which corresponds
to the discrete faithful character $\chi_0$ of the hyperbolic
structure and $m(y_0)=1$. Let $S(l,m)$ be the finite set of poles
and zeros of $l$ and $m$. Put $Y_h=\widetilde{Y}\setminus S(l,m)$
as the $X^{\prime}$ in \S 3.2. We choose the distinguished point
$t_0$ as follows. If $y_0$ is a smooth point, we take $t_0=y_0$;
if $y_0$ is a singular point, we fix a point in the pre-images of
$y_0$ in $\widetilde{Y}$ and take $t_0$ as this fixed point. This
is equivalent to fixing a branch at the singular point $y_0$.
\begin{theorem}\label{main}
(i) The closed real $1$-form $\eta(l,m)=\log{|l|}\;
d\arg{m}-\log{|m|}\; d\arg{l}$ is exact on $Y_h$;\\
(ii) For any loop $\gamma$ with initial point $t_0=\chi_0$ in
$Y_h$
\[
   \frac{1}{4
   \pi^{2}}\int_{\gamma}(\log{|m|}d\log{|l|}+\arg{l}d\arg{m})=
    \frac{k}{N},
\]
where $k$ is some integer and $N$ is the order of the symbol
$\{l,m\}$ in $K_2(\mathbb{C}(Y))$.
\end{theorem}

\begin{proof}
By (\ref{eq2.2}), we have an element $r(l,m)\in
H^{1}(Y_h,\mathbb{C}^{*})$. By Proposition \ref{prop2.1}, it is a
torsion of order $N$. By the definition of $r(l,m)$ in
(\ref{eq2.1}), we conclude that for any loop $\gamma$ in $Y_h$,
\begin{equation}\label{kt}
    \{\exp{(\frac{1}{2\pi
    i}(\int_{\gamma}\log{l}\;\frac{dm}{m}-\log{m(t_0)}\int_{\gamma}\frac{dl}{l}))}\}^N=1
\end{equation}

Write
$\displaystyle{\int_{\gamma}\log{l}\;\frac{dm}{m}-\log{m(t_0)}\int_{\gamma}\frac{dl}{l}}=Re+iIm$,
where $Re$ and $Im$ are the real and imaginary parts respectively.
(\ref{kt}) means that $\displaystyle{\exp{(\frac{N\cdot
Im}{2\pi}+\frac{N\cdot Re}{2\pi i})}=1}$. Therefore, $Im=0$ and
$\displaystyle{\frac{N\cdot Re}{2\pi i}=2 \pi i k}$, for some
integer $k$. Our result follows from the following lemma.
\end{proof}

\begin{lemma}\label{le2.1}
Denote
$\int_{\gamma}\log{l}\;\frac{dm}{m}-\log{m(t_0)}\int_{\gamma}\frac{dl}{l}=Re+iIm$
as above. Then
\[
Im=\int_{\gamma}(\log{|l|}d\arg{m}-\log{|m|}d\arg{l})=\int_{\gamma}\eta(l,m),
\]
and
\[
Re=-\int_{\gamma}(\log{|m|}d\log{|l|}+\arg{l}d\arg{m})=\int_{\gamma}\xi(l,m),
\]
where $\xi(l,m)$ depends on the branches of $\arg$ function and
$Re$ is well-defined up to $(2\pi)^{2}\mathbb{Z}$.
\end{lemma}

\begin{proof}
Let $F$ be a smooth non-zero complex-valued function, and
$F=Re(F)+iIm(F)$, where $Re(F)$ denotes its real part and $Im(F)$
its imaginary part. Then we have
\[
d\log{F}:=\frac{dF}{F}=\frac{d|F|}{|F|}+i\frac{Re(F)dIm(F)-Im(F)dRe(F)}{|F|^2}
\]

So the real part of $d\log{F}$ is $d\log{|F|}$ which is exact and
the imaginary part is denoted by $d\arg{F}$.

By a straightforward calculation, we have
\[
Im=\int_{\gamma}(\log{|l|}d\arg{m}+\arg{l}d\log{|m|})-\log{|m(t_0)|}\int_{\gamma}d\arg{l}.
\]

Integration by parts, we obtain:
\[
\int_{\gamma}\arg{l}\cdot
d\log{|m|}=\log{|m(t_0)|}\int_{\gamma}d\arg{l}-\int_{\gamma}\log{|m|}\cdot
d\arg{l}.
\]
Therefore,
\[
  Im=\int_{\gamma}(\log{|l|}d\arg{m}-\log{|m|}d\arg{l}).
\]

   For the real part $Re$, it is equal to
\[
\int_{\gamma}(\log|l|\;d\log|m|-\arg{l}\;d\arg{m})+\arg{m(t_0)}\int_{\gamma}d\arg{l}.
\]

Integration by parts, we get
\[
\int_{\gamma}\log|l|d\log|m|=-\int_{\gamma}\log{|m|}d\log{|l|}.
\]
By the choice of $t_0=\chi_{\rho_0}$, $\arg{m(t_0)}=0$. Hence the
result follows.
\end{proof}

\begin{remark}
(i) The first part of the theorem was proved in
\cite[Sect.~4.2]{CCGLS}.

(ii) The result of the second part is stronger than the one in
\cite [3.29]{Guk} where he derived that the value of the integral
is in $\mathbb{Q}$ from the quantizable Bohr-Sommerfield
condition.

(iii) The class $r(l,m)\in H^1(Y_h, \mathbb{C}^{*})$ corresponds
to a flat line bundle $L$ over $Y_h$ which is the pullback of the
universal Heisenberg line bundle on $\mathbb{C}^{*}\times
\mathbb{C}^{*}$, see \cite{Bl,Ram}. Formally,
\[
  d(\xi(l,m)+i\eta(l,m))=\frac{dl}{l}\wedge \frac{dm}{m}=0.
\]
Hence, $\frac{1}{2 \pi i}(\xi(l,m)+i\eta(l,m))$ is the $1$-form
Chern-Simons. Denote it by $CS_1(l,m)$. Then
$dCS_1(l,m)=C_1(L)=\frac{1}{2 \pi i}\frac{dl}{l}\wedge
\frac{dm}{m}=0$.
\end{remark}

By Theorem \ref{main} and Lemma \ref{le2.1}, we would like to
propose the corresponding generalized volume conjecture as the
following:\\

For a path $c: [0,1]\rightarrow Y_h$ with $c(0)=t_0$ and
$c(1)=(l,m)$, denote
\[
  U(l,m)=-\int_{c}(\log{|y|}d\log{|x|}+\arg{x}d\arg{y}).
\]
For a fixed number $a$ and $m=-\exp({i \pi a})$, we reformulate
the generalized volume conjecture as the following:
\begin{equation}\label{re-conj}
\lim_{N, k\rightarrow \infty; \frac{N}{k}=a}
\frac{\log{J_{N}(K,e^{2\pi i/k})}}{k}=\frac{1}{2 \pi}(Vol(l,m)+
i\frac{1}{2\pi}U(l,m)).
\end{equation}

\begin{remark}
By Theorem \ref{main} (ii), $\frac{1}{(2\pi)^2}U(l,m)$ is
well-defined in $\mathbb{R}/\frac{1}{N}\mathbb{Z}$. The classical
Chern-Simons invariant is well-defined in $\mathbb{R}/\mathbb{Z}$.
\end{remark}
(\ref{re-conj}) gives a $\mathbb{C}^{*}$-paramatrized version of
the volume conjecture. Using Fuglede-Kadison determinant, W. Zhang
and the first author \cite{LZ} defined an $L^2$-version twisted
Alexander polynomial which can be identified with an
$L^2$-Reidemeister torsion. By Luck and Schick's result, this
$L^2$-Alexander polynomial provides a
$\mathbb{C}^{*}$-parametrization of the hyperbolic volumes.

For other discussions on the volume conjecture, see
\cite[Sections~1.3,\;7.3]{Oh} and the related references within.

\section{$SL_2(\mathbb{C})$ Character Variety and a New Knot
Invariant} In this section, let $K$ be a hyperbolic knot in $S^3$
and $M_K$ its complement. Then $M_K$ is a hyperbolic $3$-manifold
of finite volume with boundary $\partial M_k$ a torus. Since $M_K$
is hyperbolic, there is a discrete faithful representation
$\rho_0\in R(M_K)$ corresponding to its hyperbolic structure. We
denote by $R_0$ an irreducible component of $R(M_K)$ containing
$\rho_0$. Let $X_0=t(R_0)$. By \cite{CS}, $X_0\subset X(M_K)$ is
an irreducible affine variety of dimension $1$.

\subsection{The curve of characters} In this subsection, we give some
elementary properties about the character varieties and its
component $D_0$.

Since $\partial M_K$ is a torus, we identify $R(\partial M_K)$
with the set $\{(A,B)|A,B\in SL_2(\mathbb{C}), AB=BA\}$, where
$A=\rho(\mu)$, $B=\rho(\lambda)$ for $\rho$ a representation of
$\pi_1(\partial M_K)$ in $SL_2(\mathbb{C})$. As in Section $2.3$,
$R_{D}\subset R(\partial M_K)$ is the subvariety consisting of the
representations of diagonal matrices. We have the isomorphism
$p:R_{D}\rightarrow \mathbb{C}^{*}\times \mathbb{C}^{*}$ defined
by $p(\rho)=(m,l)$ if
\begin{equation*}
\rho(\mu)=
  \left[
  \begin{matrix}
    m & 0\\
    0 & m^{-1}
  \end{matrix}
  \right], \;\;
\rho(\lambda)=
  \left[
  \begin{matrix}
    l & 0\\
    0 & l^{-1}
  \end{matrix}
  \right].
\end{equation*}

Let $R_{D}$ be identified with $\mathbb{C}^{*}\times
\mathbb{C}^{*}$ via $p$. By the proof of
\cite[Proposition~1.4.1]{CS}, $\chi\in X(\partial M_K)$ is
determined by its values on $\mu$, $\lambda$ and $\mu\lambda$.
Define a map $t:R(\partial M)\rightarrow \mathbb{C}^3$ by
$t(\rho)=(\sigma(\rho(\mu)),\sigma(\rho(\lambda)),\sigma(\rho(\mu\lambda)))$
. Then $X(\partial M_K)=t(R(\partial M_K))$. This map is the
regular surjective morphism $t:R(\partial M_K)\rightarrow
X(\partial M_K)$. That is why we use the same letter $t$. The map
$t_D$, the restriction of $t$ on $R_D=\mathbb{C}^{*}\times
\mathbb{C}^{*}$, is given explicitly by, for $(m,l)\in
\mathbb{C}^{*}\times \mathbb{C}^{*}$,
\[
t_D(m,l)=(m+m^{-1}, l+l^{-1}, ml+m^{-1}l^{-1}).
\]
It is straightforward to check that $t_D(\mathbb{C}^{*}\times
\mathbb{C}^{*})=X(\partial M_K)$ and $t_D$ is $2:1$ except at four
points $(\pm 1, \pm 1)$ where it is $1:1$.

We have the following diagram:

\[
  \begin{CD}
    R(M_K)\supset R_0   @. D_0=t^{-1}_{D}(Y_0)\subset \mathbb{C}^{*}\times \mathbb{C}^{*}=R_D\subset R(\partial M_K)\\
     @VtVV   @Vt_DVV\\
    X(M_K)\supset X_0=t(R_0) @>r>> Y_0=\overline{r(X_0)}\subset X(\partial
    M_K).
\end{CD}
\]

First we characterize the map $t_D$ and the character variety
$X(\partial M)$.
\begin{proposition}\label{f-prop}
The map $t_D:R_D=\mathbb{C}^{*}\times \mathbb{C}^{*}\rightarrow
X(\partial M_K)$ is a finite morphism.
\end{proposition}

\begin{proof}
The affine coordinate ring for $\mathbb{C}^{*}\times
\mathbb{C}^{*}$ is $\mathbb{C}[x,x^{-1},y,y^{-1}]$. Notice that
$t_D$ is surjective, we can identify the coordinate ring of
$X(\partial M_K)$ with the sub-ring
$\mathbb{C}[x+x^{-1},y+y^{-1},xy+x^{-1}y^{-1}]$ of
$\mathbb{C}[x,x^{-1},y,y^{-1}]$. Let $t=x+x^{-1}$, then $x$ and
$x^{-1}$ are roots of the equation $X^2-tX+1=0$. Hence $x$ and
$x^{-1}$ are integral over
$\mathbb{C}[x+x^{-1},y+y^{-1},xy+x^{-1}y^{-1}]$, so are $y$ and
$y^{-1}$. Now $x$, $x^{-1}$, $y$, $y^{-1}$ are integral over
$\mathbb{C}[x+x^{-1},y+y^{-1},xy+x^{-1}y^{-1}]$, it follows that
$\mathbb{C}[x,x^{-1},y,y^{-1}]$ is also integral over it.
Therefore, $t_D$ is finite.
\end{proof}

\begin{proposition}\label{torus}
The character variety $X(\partial M_K)$ is a surface in
$\mathbb{C}^3$ defined by
  \begin{equation}\label{cha-torus}
    x^2+y^2+z^2-xyz-4=0.
  \end{equation}
\end{proposition}

\begin{proof}
Let $x$, $y$, $z$ be the traces of $a=\rho(\mu)$,
$b=\rho(\lambda)$, $c=\rho(\mu\lambda)$ respectively. By the
formula (\ref{tf4}), we have
\[
  \sigma((ab)a^{-1}b^{-1})=\sigma(ab)\sigma(a^{-1}b^{-1})+\sigma(a^{-1})\sigma(a)+\sigma(b^{-1})\sigma(aba^{-1})
  -\sigma(ab)\sigma(a^{-1})\sigma(b^{-1})-\sigma(I).
\]

Since $a$ and $b$ commute, $\sigma((ab)a^{-1}b^{-1})=\sigma(I)=2$,
where $I$ is the $2\times 2$ identity matrix. This gives the
equation (\ref{cha-torus}) by (\ref{tf1}) and (\ref{tf2}). Hence,
$X(\partial M_K)$ is contained in the surface in $\mathbb{C}^3$.

On the other hand, let $(x,y,z)\in \mathbb{C}^3$ be a solution to
(\ref{cha-torus}). It is a straightforward calculation that we can
find $(m,l)\in \mathbb{C}^{*}\times \mathbb{C}^{*}$ such that
$x=m+m^{-1}$, $y=l+l^{-1}$, and $z= ml+m^{-1}l^{-1}$. Hence the
result follows.
\end{proof}

For each $\gamma \in \pi_1(M_K)$, there is a natural regular map
$I_{\gamma}: X(M_K)\rightarrow \mathbb{C}$ defined by
$I_{\gamma}(\chi)=\chi (\gamma)$. The set of functions
$I_{\gamma}$, $\gamma \in \pi_1(M_K)$ generates the affine
coordinate ring of $X(M_K)$. Moreover, by
\cite[Proposition~1.1.1]{CGLS}, for each nonzero $\gamma \in
\pi_1(\partial M_K)$, the function $I_{\gamma}$ is non-constant on
$X_0$. This implies that $r(X_0)\subset X(\partial M_K)$ has
dimension $1$. We set $Y_0=\overline{r(X_0)}$, the Zariski closure
of $r(X_0)$ in $X(\partial M_K)$. Then $Y_0$ is an irreducible
affine curve. Denote by $D_0$ the inverse image $t^{-1}_{D}(Y_0)$.

\begin{proposition}\label{d-curve}
(i) The inverse image $D_0\subset \mathbb{C}^{*}\times
\mathbb{C}^{*}$ is an affine algebraic set of dimension $1$;\\
(ii) the image of each $1$-dimensional component of $D_0$ under
$t_D$ is the whole $Y_0$;\\
(iii) $D_0$ has no $0$-dimensional components and has at most two
$1$-dimensional components.
\end{proposition}

\begin{proof}
(i). As $Y_0$ is closed, so is $D_0=t^{-1}_{D}(Y_0)$. Thus $D_0$
and $Y_0$ have the same dimension $1$ because $t_D$ is finite by
Proposition \ref{f-prop}.

(ii). Next we show that if $V$ is a $1$-dimensional irreducible
component of $D_0$, then $t_D(V)=Y_0$. In fact, since a finite
morphism maps closed sets to closed sets, $t_D(V)\subset Y_0$ is
closed, irreducible and its dimension is $1$. Hence $t_D(V)=Y_0$
due to the irreducibility of $Y_0$.

(iii). Suppose that $D_0$ has three distinct $1$-dimensional
components $V_i$, $1\leq i \leq 3$. Since the $V_i \cap V_j$,
$i\ne j$ are empty or finite sets, we can choose $y_0\in Y_0$ such
that $y_0\notin t_D(V_i \cap V_j)$ for $\forall i\ne j$. We see
that $t_D(V_i)=Y_0$, hence $t_D^{-1}(y_0)$ has three elements.
This contradicts that $t_D^{-1}(y)$ has at most two elements for
any $y\in Y_0$. Therefore, $D_0$ has at most two $1$-dimensional
irreducible components.

Now we show that it has no $0$-dimensional components. Consider
the morphism $\sigma: \mathbb{C}^{*}\times
\mathbb{C}^{*}\rightarrow \mathbb{C}^{*}\times \mathbb{C}^{*}$
defined by $\sigma(m,l)=(m^{-1},l^{-1})$. It is an involution and
hence an automorphism of $\mathbb{C}^{*}\times \mathbb{C}^{*}$
with $\sigma(D_0)=D_0$. So it is an automorphism of $D_0$ when
restricting on it. Therefore it maps $1$-dimensional component to
$1$-dimensional component. Now for any $y\in Y_0$, $\sigma$
permutes elements of $t_D^{-1}(y)$ and at least one element of
$t_D^{-1}(y)$ is contained in some $1$-dimensional component since
$t_D$ maps every $1$-dimensional component onto $Y_0$. Therefore,
no $0$-dimensional component for $D_0$.
\end{proof}

\begin{remark}
For the character variety $X(M_K)$, by
\cite[Proposition~2.4]{CCGLS}, there is no zero-dimensional
components of $X(M_K)$. Proposition \ref{d-curve} part (iii) shows
that $D_0$ in $R_D\subset R(\partial M_K)$ has no $0$-dimensional
components.
\end{remark}

In summary, we have the following two cases:\\
(I) $D_0$ itself is an irreducible affine curve;\\
(II) $D_0=V_1\cup V_2$, where $V_i$ are two $1$-dimensional
irreducible components. Each is an irreducible affine curve with
$t_D(V_i)=Y_0$, $i=1,2$. Moreover, they are isomorphic to each
other under the involution $\sigma$ by Proposition \ref{d-curve}
(iii) . We also have $t_D: V_i\rightarrow Y_0$ is a one-to-one and
onto regular map.

\subsection{Invariant and Dehn Surgery of Hyperbolic Knots}
In this subsection, we define an algebro-geometric invariant of
the $3$-manifolds $K(p/q)$ resulting from the $(p,q)$ Dehn surgery
along a hyperbolic knot complement in $S^3$. A Casson type
$SL_2(\mathbb{C})$ invariant for $K(p/q)$ is also established.

Let $A_0(m,l)$ be the defining equation of the closure of the
affine curve $D_0$ in $\mathbb{C}\times \mathbb{C}$. We also
require that it has no repeated factors. It is a factor of the
$A$-polynomial of the knot $K$ defined in \cite{CCGLS}, and
$A_0(m,l)^d$ is the $A$-polynomial of $X_0$ defined in
\cite[Page~109]{BZ}, where $d$ is the degree of regular map
$r:X_0\rightarrow Y_0$.

Denote by $\widetilde{X_0}$ (resp. $\widetilde{Y_0}$) a smooth
projective model of the affine curve $X_0$ (resp. $Y_0$). The
restriction morphism $r:X_0\rightarrow Y_0$ induces a regular map
$\widetilde{r}: \widetilde{X_0}\rightarrow \widetilde{Y_0}$.

\begin{lemma}\label{le1}
The regular map $\widetilde{r}$ is an isomorphism.
\end{lemma}

\begin{proof}
By our assumption, $M_K$ is the complement of a knot in $S^3$,
hence $H^1(M_K;\mathbb{Z}_2)=\mathbb{Z}_2$. By
\cite[Corollary~3.2]{Dun}, $r$ is a birational isomorphism onto
$r(X_0)$. Since $\widetilde{r}$ is induced by $r$, the result
follows.
\end{proof}

Let $\gamma=p\mu+q\lambda\in H_1(\partial M_K;\mathbb{Z})$ be a
non-zero primitive element with $p,q$ coprime. Define a regular
function $f_{\gamma}=I_{\gamma}^2-4$ on $X_0$. Since $I_{\gamma}$
is nonconstant, so is $f_{\gamma}$. It is also a meromorphic
function on $\widetilde{X_0}$ or equivalently, a non-constant
holomorphic function from $\widetilde{X_0}$ to $\mathbb{CP}^1$ and
we denote it again by $f_{\gamma}$ .

Since $\gamma\in H_1(\partial M_K;\mathbb{Z})$, we can think of
$I_{\gamma}$ as a regular function on $Y_0$ in $X(\partial M)$.
Define on $Y_0$ the function $f_{\gamma}^{\prime}=I_{\gamma}^2-4$.
Similarly, it is a non-constant regular function on $Y_0$, and
hence a non-constant holomorphic function from $\widetilde{Y_0}$
to $\mathbb{CP}^1$, denoted also by $f_{\gamma}^{\prime}$. Then by
the definition, we have
\[
f_{\gamma}=f_{\gamma}^{\prime}\circ \widetilde{r}.
\]
In particular,
$\text{deg}f_{\gamma}=\text{deg}f_{\gamma}^{\prime}$ by Lemma
\ref{le1}.

We denote by $Z_{\gamma}$ the set of zeroes of the function
$f_{\gamma}$ on $X_0$. If $\chi\in Z_{\gamma}$, then there exists
a representation $\rho\in R_0$ such that its character
$\chi_\rho=\chi$.

\begin{lemma}\label{le2}
Suppose $\chi_\rho\in Z_{\gamma}$. Then either $\rho(\gamma)=\pm
\text{I}$ or $\rho(\gamma)\ne \pm I$ and $\sigma(\rho(\alpha))=\pm
2$ for all $\alpha\in \pi_1(\partial M_K)$.
\end{lemma}

\begin{proof}
Note that $f_{\gamma}(\chi_\rho)=0$ is equivalent to that the
trace of $\rho(\gamma)$ is $\pm 2$. If $\rho(\gamma)\ne \pm I$,
then $\rho(\gamma)$ is a parabolic element in $SL_2(\mathbb{C})$.
Since $\rho(\gamma)$ and $\rho(\alpha)$ commute for all $\alpha\in
\pi_1(\partial M_K)$, $\sigma(\rho(\alpha))=\pm 2$.
\end{proof}

Suppose that $\rho\in R_0\subset R(M_K)$ is irreducible. Then its
character $\chi_{\rho}$ is contained in the component $X_0$.
Assume that $\rho(\mu)$ and $\rho(\lambda)$ are parabolic. Up to
conjugation, we have
\begin{equation*}
\rho(\mu)=\pm
  \left[
  \begin{matrix}
    1 & 1\\
    0 & 1
  \end{matrix}
  \right]; \;\;
\rho(\lambda)=\pm
  \left[
  \begin{matrix}
    1 & t(\rho)\\
    0 & 1
  \end{matrix}
  \right],
\end{equation*}

where $t(\rho)$ is a complex number. We have the following:
\begin{conjecture}
Let $\rho\in R_0$ be an irreducible
$SL_2(\mathbb{C})$-representation of a hyperbolic knot in $S^3$.
If $\rho(\mu)$ and $\rho(\lambda)$ are parabolic, then
$t(\rho)\notin \mathbb{Q}$.
\end{conjecture}

\begin{remark}
(1) If $\rho=\rho_0$ is the discrete faithful representation of
the hyperbolic structure, then we know $\rho(\mu)$ and
$\rho(\lambda)$ are parabolic and $t(\rho_0)$ is called the
\emph{cusp constant} and the \emph{cusp polynomial} is the minimum
polynomial for $t(\rho_0)$ over $\mathbb{Q}$. Moreover $\{1,
t(\rho_0)\}$ generates a lattice of $\mathbb{C}$. Therefore,
$t(\rho_0)\notin \mathbb{R}$. For more detail, see
\cite[Section~6]{CL}. So $t(\rho)$ can be thought of as the
generalization of the cusp constant.

(2) The conjecture can not be extended to non-hyperbolic knots.
There is an example in \cite[Section~5]{Ril1} of an irreducible
parabolic representation $\rho$ of alternating torus knot such
that $t(\rho)\in \mathbb{Z}$.
\end{remark}

The conjecture is true for the figure-eight knot. In fact, for the
figure-eight knot, by Proposition \ref{8cha-1}, if $\rho$
satisfies the condition of the conjecture, then it must be the
discrete faithful representation of the hyperbolic structure. By
the remark, $t(\rho)\notin \mathbb{R}$.

\begin{proposition}\label{kin}
Suppose the above conjecture is true. Let $\gamma=p\mu+q\lambda\in
H_1(\partial M_K;\mathbb{Z})$ with $p,q$ non-zero coprime
integers. Let $\chi_{\rho}\in Z_{\gamma}$ be the character of an
irreducible representation $\rho$. Then $\rho(\gamma)=\pm I$ if
and only if the trace of $\rho(\mu)$ is not equal to $\pm 2$.
\end{proposition}

\begin{proof}
If the trace of $\rho(\mu)$ is not equal to $\pm 2$, then
$\rho(\mu)$ is not parabolic or $\pm I$. Since $\rho(\mu)$ and
$\rho(\gamma)$ commute, $\rho(\gamma)$ is not parabolic. Since
$\chi_{\rho}\in Z_{\gamma}$, $\sigma(\rho(\gamma))=\pm 2$. Thus,
we obtain $\rho(\gamma)=\pm I$.\\

Next suppose that the trace of $\rho(\mu)$ is equal to $\pm 2$.
We show that $\rho(\gamma)\ne \pm I$. There are two cases.\\
Case I. $\rho(\mu)=\pm I$. We know that the $1/0$ Dehn surgery
produces $S^3$. Hence $\rho$ induced a representation of
$\pi_1(S^3)$ in $PSL_2(\mathbb{C})$. It is trivial. Therefore, the
image of $\rho$ in $SL_2(\mathbb{C})$ is contained in $\{\pm I\}$.
So $\rho$ is reducible. This contradiction shows that Case I can
not happen.\\
Case II. $\rho(\mu)$ is parabolic. Since $\rho(\mu)$ and
$\rho(\lambda)$ commute, up to conjugation, we can assume that
\begin{equation*}
\rho(\mu)=\pm
  \left[
  \begin{matrix}
    1 & 1\\
    0 & 1
  \end{matrix}
  \right], \;\;
\rho(\lambda)=\pm
  \left[
  \begin{matrix}
    1 & t(\rho)\\
    0 & 1
  \end{matrix}
  \right].
\end{equation*}

Now $\gamma=p\mu+q\lambda$, so
$\rho(\gamma)=\rho(\mu)^p\rho(\lambda)^q$. If $t(\rho)=0$,
$\rho(\gamma)\ne \pm I$ unless $p=0$. If $t(\rho)\ne 0$, by the
assumption $t(\rho)\notin \mathbb{Q}$, then $\rho(\gamma)\ne \pm
I$.
\end{proof}

\begin{proposition}\label{redu}
Suppose $\chi\in Z_{\gamma}$ is the character of a reducible
representation $\rho$. Then $\chi(\mu)\ne \pm 2$ and
$\rho(\gamma)=\pm I$.
\end{proposition}

\begin{proof}
Let $m$ be an eigenvalue of $\rho(\mu)$. By
\cite[Proposition~6.2]{CCGLS}, $m^2$ must be a root of the
Alexander polynomial $\Delta(t)$ of the knot. It is well-known
that $\Delta(1)\ne 0$, hence $m\ne \pm 1$ and $\chi(\mu)\ne \pm
2$. This also means that $\rho(\mu)$ is not parabolic or $\pm I$.
Since $\rho(\mu)$ and $\rho(\gamma)$ commute and the trace of
$\rho(\gamma)$ is $\pm 2$, we must have $\rho(\gamma)=\pm I$.
\end{proof}

\begin{remark}
Compare with Proposition \ref{kin}, for the reducible characters,
it is much simpler and there is no need of the conjecture.
\end{remark}

Let $E(p,q)$ be the reducible curve $m^pl^q=\pm 1$ in
$\mathbb{C}^{*}\times \mathbb{C}^{*}$ for $p,q$ coprime integers.
Then the image $t_D(E(p,q))$ is a curve in $X(\partial M_K)$. We
know $r(X_0)$ is an irreducible curve in $X(\partial M_K)$. They
do not have common irreducible component because the traces of
characters of $X_0$ are not constant. Hence $t_D(E(p,q))\cap
r(X_0)$ is finite. The set $t_D(E(p,q))\cap r(X_0)$ consists of
possible characters in $X_0$ which can also be the characters of
$K(p/q)$, where $K(p/q)$ denotes the closed $3$-manifold obtained
from $M_K$ by the Dehn surgery along the simple closed curve of
$\partial M_K$ which represents the class $\gamma$ in
$H_1(\partial M_K;\mathbb{Z})$. The following definition should be
thought of as the \emph{algebro-geometric} invariant for the
$(p,q)$ Dehn surgery of $M_K$.

\begin{definition}
\[
  b(p,q)=\sum_{\chi\in t_D(E(p,q))\cap r(X_0)}n_{\chi},
\]
where $n_{\chi}$ is the intersection multiplicity at $\chi$.
\end{definition}

\begin{theorem}\label{b-inv}
The integer $b(p,q)$ is a well-defined invariant of the
$3$-manifold $K(p/q)$ resulting from the Dehn filling on the
hyperbolic knot complement $M_K$. It is always positive.
\end{theorem}

\begin{proof}
As mentioned above, the set $t_D(E(p,q))\cap r(X_0)$ is finite and
hence $b(p,q)<\infty$. On the other hand, because
$\chi_{\rho_0}(\mu)=\pm 2$, $\chi_{\rho_0}$ is always contained in
this set, where $\rho_0\in R_0$ is the discrete faithful
representation of the hyperbolic metric. Therefore, $b(p,q)>0$.

The intersection between $t_D(E(p,q))$ and $r(X_0)$ is taking
place in the surface $X(\partial M_K)$ described in Proposition
\ref{torus}. If two hyperbolic knots $K_1$ and $K_2$ are
homeomorphic, then they have the isomorphic fundamental groups of
their complements in $S^3$, hence they have isomorphic $X_0$.
Therefore, they have the same $b(p,q)$. Thus $b(p,q)$ is an
invariant of $K(p/q)$ depending only on the hyperbolic knot $K$
and the Dehn surgery coefficient $p/q$.
\end{proof}

Now set $S(p,q)=\{\chi\in t_D(E(p,q))\cap r(X_0)|\chi(\mu)\ne \pm
2\}\subset r(X_0)\subset X(\partial M)$.

\begin{proposition}\label{count}
Suppose the conjecture 1 is true. Then $S(p,q)$ is exactly the set
of characters in $X_0$ which are also the characters of $K(p/q)$.
\end{proposition}

\begin{proof}
This follows from the definition of $S(p,q)$ and Proposition
\ref{kin}.
\end{proof}

\begin{definition}\label{lambda}
$$\lambda (p,q)=\sum_{\chi\in S(p,q)} n_{\chi},$$
where $n_{\chi}$ is the intersection multiplicity at $\chi$.
\end{definition}

By Proposition \ref{redu}, the set $S(p,q)$ contains all possible
reducible characters. Hence the number $\lambda (p,q)$ counts both
irreducible and reducible characters of $K(p/q)$. From the
mathematical physics point of view in \cite{Guk}, it is important
to count both irreducible and reducible representations in
$SL_2(\mathbb{C})$ case. It is easy to count the abelian and
non-abelian reducible characters, so there is a computable way to
count the irreducible characters for $K(p/q)$.

\begin{theorem}\label{c-inv}
Assume the conjecture 1 is true. The quantity $\lambda (p,q)$ is a
well-defined \emph{algebro-geometric} $SL_2(\mathbb{C})$ Casson
type invariant of $K(p/q)$.
\end{theorem}

\begin{proof}
It follows from Theorem \ref{b-inv} and Proposition \ref{count}.
\end{proof}

Note that by definition, $\lambda (p,q)\leq b(p,q)$ for any
coprime $p$, $q$ and a hyperbolic knot in $S^3$.

The invariant $\lambda (p,q)$ is for the hyperbolic knot and its
$(p,q)$ Dehn surgery. An $SL_2(\mathbb{C})$ knot invariant
obtained from the character variety of $1$-dimensional components
is given by the first author in \cite{Li}. The construction in
\cite{Li} was purely topological by choosing generic smooth
perturbations and generic values of $\chi(\mu)$. The topological
definition of the Casson-type invariant is not easy to calculate.
Our invariant $\lambda (p,q)$, defined via the intersection
multiplicity, is easier or at least very explicit for computation.
It can also be interpreted as the intersection cycles of the
appropriate cohomology classes of $X(\partial M)$.

\begin{proposition}\label{prop-liv}
If $\lambda (p,q)>0$, then $\pi_1(K(p/q))$ is non-cyclic.
\end{proposition}

\begin{proof}
(1) If $\lambda (p,q)>0$, then there exists $\rho\in R_0$ such
that its character $\chi_{\rho}\in S(p,q)$. Hence
$F_{\mu}(\chi_{\rho})\ne 0$, and $F_{\gamma}(\chi_{\rho})=0$. By
\cite[Proposition~1.5.2]{CGLS}, there exists a representation from
$\pi_1(K(p/q))$ to $PSL_2(\mathbb{C})$ with non-cyclic image. The
result follows.
\end{proof}

\begin{remark}
(1) Our invariant $\lambda (p,q)$ is defined as the algebraic
intersection multiplicity in $X(\partial M)$ from the $X_0$
component. This is different from the definition in \cite{Cu1}
where the number $\lambda(K_{p/q})$ is defined over all components
of $X(M)$ and the intersection is taken in different space.\\
(2) For hyperbolic knot $K$, $K(p/q)$ may not be NSL manifold. For
the component $X_0$, the intersection in \cite{Cu1} is $r(X_0)\cap
(t_D(t_D^{-1}(\overline{r(X_0)})\cap \{m^pl^q=1\}))$. This is
different from our Definition \ref{lambda}. Moreover, the
intersection multiplicity of an intersection point in $X(\partial
M_K)$ is counted by its multiplicity in $X(K(p/q))$ via an
appropriate Heegaard splitting of $K(p/q)$ in \cite{Cu1}. \\
(3) It would be interesting to prove that the two character
varieties of Heegaard handle-bodies are smooth and their
intersection is always proper. In \cite[Section~4]{FM}, Fulton and
MacPherson only sketched an argument that if the two smooth
subvarieties intersect properly, then the topological intersection
multiplicity agrees with the algebraic intersection multiplicity .
\end{remark}

In \cite[Section~1.4]{CGLS}, a norm $|.|$ is defined on the real
vector space $H_1(\partial M_K,\mathbb{R})$ with the property that
$|\gamma|=\text{deg}f_{\gamma}$ for any $\gamma\in H_1(\partial
M_K,\mathbb{Z})$. This norm is called the Culler-Shalen norm. In
particular, $|\gamma|=\text{deg}f_{\gamma}^{\prime}$.

Let $\pi: \widetilde{Y_0}\rightarrow Y_0$ be the birational
isomorphism. Note that $\pi$ is well-defined only on a Zariski
dense subset of $\widetilde{Y_0}$ and is surjective. A point of
the set $I=\widetilde{Y_0}\setminus  \pi^{-1}(Y_0)$ is called an
ideal point. Denote by $Z_{\gamma}^{\prime}$ the set of zeroes of
the meromorphic function
$f_{\gamma}^{\prime}:\widetilde{Y_0}\rightarrow \mathbb{C}$. Then
\[
  \text{deg}f_{\gamma}^{\prime}=\sum_{y\in
  Z_{\gamma}^{\prime}}v_y,
\]
where $v_y$ is the order of vanishing of $f_{\gamma}^{\prime}$ at
$y$.\\

Set $Z_1=\{y\in Z_{\gamma}^{\prime}|\pi(y)\in S(p,q)\}$ and
$I(p,q)=Z_{\gamma}^{\prime}\cap I$. Now we have another quantity:
\begin{equation}\label{v-num}
  \widehat{\lambda}(p,q)=\sum_{y\in Z_1}v_y.
\end{equation}
The natural question is to find the relationship between
$\widehat{\lambda}(p,q)$ and $\lambda(p,q)$ of the Definition
\ref{lambda}. It seems that there is no easy answer to this
question. See the remark below. Nevertheless, we have the
following:
\begin{proposition}\label{deg-ine}
Assume that for every $\chi\in S(p,q)$, $t_D(E(p,q))$ intersects
$r(X_0)$ transversely at $\chi$. Then\\
(i) $\lambda(p,q)\leq \widehat{\lambda}(p,q)$,\\
(ii) $\lambda(p,q)+\widehat{I}(p,q)\leq \text{deg}f_{\gamma}$,
where $\widehat{I}(p,q)=\sum_{x\in I(p,q)}v_x$.
\end{proposition}

\begin{proof}
If the intersection is transverse, then intersection multiplicity
$n_{\chi}=1$. The order of vanishing is at least one,  hence (i)
holds.

For (ii), note that $Z_1$ and $I(p,q)$ are subsets of
$Z_{\gamma}^{\prime}$. Thus, $\lambda(p,q)+\widehat{I}(p,q)\leq
\text{deg}f_{\gamma}^{\prime}=\text{deg}f_{\gamma}$.
\end{proof}

\begin{remark}
If the intersection is not transversal, for instance, some $\chi$
is a singular point of $r(X_0)$, then we do not know how to
compare them. One difficulty is that when $\chi$ is a singular
point, there is NO well-defined notion of the order of vanishing
of $f_{\gamma}^{\prime}$ at $\chi$ in $Y_0$. Moreover,
$\pi^{-1}(\chi)$ has more than one element in $\widetilde{Y_0}$
and each of them is a zero of $f_{\gamma}^{\prime}$.
\end{remark}

\section{An Example: The Figure Eight Knot} Throughout this
section, we shall denote by $M$ the complement of the figure-eight
knot in $S^3$. we compute its character variety and give some
applications.

It is well-known that $\pi_1(M)$ is given by two generators and
one relation:
\begin{equation}\label{8gp}
  \pi_1(M)=<\alpha,\beta|R(\alpha,\beta)>,\;
  R(\alpha,\beta)=\beta^{-1}\alpha^{-1}\beta\alpha\beta^{-1}\alpha\beta\alpha^{-1}\beta^{-1}\alpha,
\end{equation}

where $\alpha$, $\beta$ are meridians, and they are conjugate to
each other $\beta=\delta\alpha\delta^{-1}$ with
$\delta=\alpha^{-1}\beta\alpha\beta^{-1}$.

Set $\tau=\alpha\beta^{-1}\alpha^{-1}\beta$, we have a peripheral
subgroup $\pi_1(\partial M)$
\begin{equation}\label{8bgp}
 \pi_1(\partial M)=<\alpha,\lambda>; \;
 \lambda=\tau^{-1}\delta=\beta^{-1}\alpha\beta\alpha^{-1}\alpha^{-1}\beta\alpha\beta^{-1},
\end{equation}

where $\alpha$ is the meridian and $\lambda$ is the longitude.
Note in this section, we use $\alpha$ for the meridian instead of
$\mu$.

Let us consider its character variety $X(M)$. We use (\ref{8gp})
for the presentation of $\pi_1(M)$. For a
$SL_2(\mathbb{C})$-representation $\rho\in R(M)$, its character
$\chi$ is determined by $\chi(\alpha)$, $\chi(\beta)$, and
$\chi(\alpha\beta)$. We have the morphism $t:R(M)\rightarrow
\mathbb{C}^3$, $t(\rho)=(\chi(\alpha),
\chi(\beta),\chi(\alpha\beta))$ and $X(M)$ is the image $t(R(M))$.
Let $(x,y,z)$ be the affine coordinate for $\mathbb{C}^3$.

\begin{proposition}\label{8cha}
The affine variety $X(M)\subset \mathbb{C}^3$ is defined by the
following equations:
\begin{equation}\label{8d1}
  x=y,
\end{equation}
\begin{equation}\label{8d2}
 (x^2-z-2)(z^2-(1+x^2)z+2x^2-1)=0.
\end{equation}
In particular, we can identify $X(M)$ with the affine plane curve
$\{(x,z)\in \mathbb{C}^2|(x^2-z-2)(z^2-(1+x^2)z+2x^2-1)=0\}$. It
has two irreducible components.
\end{proposition}

\begin{proof}
Since $\alpha$ is conjugate to $\beta$,
$x=\chi(\alpha)=\chi(\beta)=y$. Hence (\ref{8d1}) follows. For the
second equation, by \cite[Theorem~1]{Wh}, the factor $x^2-z-2$
corresponds to characters of abelian representations, and the
other one corresponds to the characters of non-abelian
representations.
\end{proof}

By the preceding proposition, the component defined by the
equation $z^2-(1+x^2)z+2x^2-1=0$ consists of the characters of
non-abelian representations, in particular, it contains the
discrete faithful representations of the complete hyperbolic
metric of $M$. We denote this component $X_0$. Hence we have :

\begin{corollary}\label{8-cor1}
The component $X_0$ is an irreducible smooth curve in
$\mathbb{C}^2$ with the defining equation:
\begin{equation}\label{8d3}
   z^2-(1+x^2)z+2x^2-1=0.
\end{equation}
\end{corollary}

\begin{proof}
We only need to check that (\ref{8d3}) defines a smooth curve. Let
$f(x,z)=z^2-(1+x^2)z+2x^2-1$. Then $\frac{\partial f}{\partial
x}=-2xz+4x$ and $\frac{\partial f}{\partial z}=2z-1-x^2$. It is
straightforward to check that their is no common solution to the
equations $\frac{\partial f}{\partial x}=\frac{\partial
f}{\partial z}=f(x,z)=0$. Hence, the curve is smooth.
\end{proof}

\begin{corollary}\label{8-cor2}
There are exactly two reducible characters on $X_0$ and they
correspond to the points $(\pm \sqrt{5},3)$.
\end{corollary}

\begin{proof}
Let $\rho$ be a reducible representation whose character lies on
$X_0$. Let $m$ be the eigenvalue of $\rho(\alpha)$. By
\cite[Proposition~6.2]{CCGLS}, $m^2$ must be a root of the
Alexander polynomial $\Delta(t)=t^2-3t+1$ of the figure-eight
knot. Thus,
\[
  m=\pm \sqrt{\frac{3\pm \sqrt{5}}{2}}=\frac{\pm \sqrt{5}\pm
  1}{2}; \, \text{and}\; x=\sigma(\rho(\alpha))=m+m^{-1}=\pm \sqrt{5}.
\]

Plug in (\ref{8d3}), we get $z=3$ with multiplicity two. Now it is
easy to check that $(\pm \sqrt{5},3)$ are exactly the intersection
points of $X_0$ and the component $x^2-z-2=0$. Since the latter
component consists of abelian characters, the result follows.
\end{proof}

\begin{remark}
Two inequivalent reducible representations may have the same
character. The points $(\pm \sqrt{5},3)$ above are such examples.
They are the characters of both some abelian representation and
non-abelian reducible representation which are clearly not
equivalent.
\end{remark}

\begin{proposition}\label{8cha-1}
Suppose that $\chi\in X_0$ and $\chi(\alpha)=\pm 2$. Then $\chi$
is the character of a discrete faithful representation.
\end{proposition}

\begin{proof}
Plug $x=\pm 2$ in (\ref{8d3}), we get $z^2-5z+7=0$ and hence
$z=\frac{5 \pm \sqrt{-3}}{2}$. The following representation
$\theta:\pi_1(M)\rightarrow PSL_2(\mathbb{C})$ was given in
\cite{Ril2}:
\begin{equation*}
\theta(\alpha)=\pm
  \left[
  \begin{matrix}
    1 & 1\\
    0 & 1
  \end{matrix}
  \right], \;\;
\theta(\beta)=\pm
  \left[
  \begin{matrix}
    1 & 0\\
    -\omega & 1
  \end{matrix}
  \right];
\end{equation*}
where $\omega=\frac{-1+\sqrt{-3}}{2}$ is a primitive cube root of
unity. By \cite[Theorem~1]{Ril2}, $\theta$ is an isomorphism.
Thus, $\theta$ is a discrete faithful representation of
$\pi_1(M)$. Let $\overline{\theta}$ be the complex conjugation of
$\theta$. It is straightforward to check that the four points
$(\pm 2, \frac{5 \pm \sqrt{-3}}{2})$ correspond to the
$SL_2(\mathbb{C})$ characters of the lifts of $\theta$ and
$\overline{\theta}$. The result follows.
\end{proof}

Denote by $\mathbb{CP}^2$ the complex projective plane. We use
$[X:Y:Z]$ to represent its homogenous coordinates. We will
identify $\mathbb{C}^2$ with the open subset $\{[x:1:z]|(x,z)\in
\mathbb{C}^2\}$. Note this is different from the standard
notation.

Let $\widetilde{X_0}$ be a smooth projective model of $X_0$. We
have the following explicit description of $\widetilde{X_0}$ in
$\mathbb{CP}^2$.
\begin{proposition}\label{8-smo}
$\widetilde{X_0}$ is an elliptic curve and defined by the equation
in $\mathbb{CP}^2$:
\begin{equation}\label{8sm}
   YZ^2-Y^2Z-X^2Z+2X^2Y-Y^3=0.
\end{equation}
\end{proposition}

\begin{proof}
Equation (\ref{8d3}) is the defining equation of $X_0$. Substitute
$x=\frac{X}{Y}$, $z=\frac{Z}{Y}$, and we get (\ref{8sm}). That is,
it is the projective closure of $X_0$ in $\mathbb{CP}^2$. It
suffices to show that it is smooth. Let $F(X,Y,Z)$ be the
left-hand side of (\ref{8sm}). It is elementary to check that
except $(0,0,0)$, there is no common solution to the equations
$\frac{\partial F}{\partial X}=\frac{\partial F}{\partial
Y}=\frac{\partial F}{\partial Z}=0$. Hence, it is smooth. Since
(\ref{8sm}) is a cubic equation in $\mathbb{CP}^2$,
$\widetilde{X_0}$ has genus $1$.
\end{proof}

\begin{corollary}\label{8-ide}
$X_0$ has exactly two ideal points $[1:0:0]$ and $[0:0:1]$.
\end{corollary}

\begin{proof}
The ideal points correspond to $Y=0$. This gives $X^2Z=0$. Thus,
$X=0$ or $Z=0$.
\end{proof}

 Next, we explicitly construct the restriction map
$r:X_0\rightarrow X(\partial M)$ induced by the inclusion
$i:\pi_1(\partial M)\rightarrow \pi_1(M)$.

\begin{theorem}\label{8res}
The map $r:X_0\rightarrow X(\partial M)\in \mathbb{C}^3$ is given
by the formulas:
\[
  r(x,z)=(x,F(x,z),G(x,z));
\]
where $F(x,z)=x^4-5x^2+2$ is the trace of the longitude $\lambda$,
and $G(x,z)=(4x-x^3)z+(x^5-4x^3-x)$ is the trace of $\alpha
\lambda$.
\end{theorem}

We will postpone the proof of Theorem \ref{8res} to the end of
this section. Instead, we first discuss its applications. The
reason is that the proof itself is elementary and lengthy, and
probably is known to the experts. However, we can not find it in
the literature. So we include here for completeness.

Denote by $x_0=[0:0:1]$ and $x_1=[1:0:0]$. By Corollary
\ref{8-ide}, they are the only ideal points of $X_0$.

\begin{lemma}\label{l-norm}
(1) Both $x_0$ and $x_1$ are poles of the meromorphic functions
$f_{\alpha}$ and $f_{\lambda}$. Their orders are:
$v_{x_0}(f_{\alpha})=v_{x_1}(f_{\alpha})=2$, and
$v_{x_0}(f_{\lambda})=v_{x_1}(f_{\lambda})=8$. \\
(2) The Culler-Shalen norm $|(1,0)|=4$ and $|(0,1)|=16$.
\end{lemma}

\begin{proof}
(1) By Theorem \ref{8res}, for $\chi=(x,z)\in X_0$,
$I_{\alpha}(x,z)=x$ and $I_{\lambda}(x,z)=F(x,z)=x^4-5x^2+2$.
Since $f_{\alpha}=I_{\alpha}^2-4$ and
$f_{\lambda}=I_{\lambda}^2-4$, it is sufficient to show that
$v_{x_0}(I_{\alpha})=v_{x_1}(I_{\alpha})=1$, and
$v_{x_0}(I_{\lambda})=v_{x_1}(I_{\lambda})=4$.

On $\widetilde{X_0}$, substitute $x=\frac{X}{Y}$, $z=\frac{Z}{Y}$,
we get for $[X:Y:Z]\in \widetilde{X_0}$,
\[
   I_{\alpha}([X:Y:Z])=\frac{X}{Y},\, \text{and}\,\,
   I_{\lambda}([X:Y:Z])=\frac{X^4-5X^2Y^2+2Y^4}{Y^4}.
\]
First, we consider $x_0=[0:0:1]$. Let $U_3=\{[X:Y:Z]|Z\ne 0\}$.
Then it is an open subset containing $x_0$. $U_3$ is identified
with $\mathbb{C}^2$ via $x=\frac{X}{Z}$ and $y=\frac{Y}{Z}$, where
$(x,y)$ are affine coordinates of $\mathbb{C}^2$. Dividing by
$Z^3$  the both sides of the equation (\ref{8sm}) of
$\widetilde{X_0}$ and substituting $x=\frac{X}{Z}$ and
$y=\frac{Y}{Z}$, we get the affine equation in $U_3=\mathbb{C}^2$:
\[
  g(x,y)=y-y^2-x^2+2x^2y-y^3=0.
\]
Under this identification, $x_0$ is the origin $(0,0)$,
$I_{\alpha}(x,y)=\frac{x}{y}$ and
$I_{\lambda}(x,y)=(\frac{x}{y})^4-5(\frac{x}{y})^2+2$. Since
$\frac{\partial g}{\partial y}=1-2y+2x^2-3y^2$ and $\frac{\partial
g}{\partial y}(0,0)\ne 0$, the function $x$ is a local parameter
of the local ring of regular functions at $(0,0)$. Solve
$g(x,y)=0$, we get
\[
  I_{\alpha}(x,y)=\frac{x}{y}=x^{-1}u(x,y),\, \text{where}\,
  u(x,y)=\frac{y^2+y-1}{2y-1},\, u(0,0)\ne 0;
\]
and
\[
  I_{\lambda}(x,y)=x^{-4}w(x,y),\, \text{where}\,
  w(x,y)=u^4-5x^2u^2+2x^4
  ,\, w(0,0)\ne 0.
\]
Therefore, $x_0$ is a pole of $I_{\alpha}$ of order $1$, and it is
a pole of $I_{\lambda}$ of order $4$.

For $x_1=[1:0:0]$. Let $U_1=\{[X:Y:Z]|X\ne 0\}$. $U_1$ is
identified with $\mathbb{C}^2$ via $y=\frac{Y}{X}$ and
$z=\frac{Z}{X}$, where $(y,z)$ are affine coordinates of
$\mathbb{C}^2$. Similarly, we obtain the affine equation in
$U_1=\mathbb{C}^2$:
\[
  h(x,y)=yz^2-y^2z-z+2y-y^3=0.
\]
Now $x_0$ is the origin $(0,0)$, $I_{\alpha}(x,y)=y^{-1}$ and
$I_{\lambda}(x,y)=y^{-4}-5y^{-1}+2$. Since $\frac{\partial
h}{\partial z}(0,0)\ne 0$, the function $y$ is a local parameter
of the local ring of regular functions at $(0,0)$. Thus, $x_0$ is
a pole of $I_{\alpha}$ of order $1$, and it is a pole of
$I_{\lambda}$ of order $4$.

(2) By definition,
$|(1,0)|=\text{deg}f_{\alpha}=v_{x_0}(f_{\alpha})+v_{x_1}(f_{\alpha})$,
by (1),  $|(1,0)|=2+2=4$. Similarly, $|(0,1)|=8+8=16$.
\end{proof}

Use this lemma, we can compute the Culler-Shalen norm $|(p,q)|$
for any $\gamma=p\alpha+q\lambda\in H_1(\partial M,\mathbb{Z})$.

\begin{proposition}\label{p-8norm}
For each $\gamma=p\alpha+q\lambda\in H_1(\partial M,\mathbb{Z})$,
the Culler-Shalen norm $|\gamma|=|(p,q)|=2(|p+4q|+|p-4q|)$.
\end{proposition}

\begin{proof}
By Lemma \cite[~1.4.1, 1.4.2]{CGLS}, for each ideal point $x$,
there is a homomorphism $\phi_x:H_1(\partial
M,\mathbb{Z})\rightarrow \mathbb{Z}$, such that, for each
$\gamma=p\alpha+q\lambda$,
\[
  v_x(f_{\gamma})=|\phi_x(\gamma)|,\, \text{and}\, |\gamma|=\sum_{x:\,\text{ideal
  point}}|\phi_x(\gamma)|;
\]
where $v_x(f_{\gamma})$ denotes the order of pole of $f_{\gamma}$
at $x$. For our case, we have two ideal points $x_0$ and $x_1$. By
Lemma \ref{l-norm}, we get $|\phi_{x_0}(\alpha)|=2$,
$|\phi_{x_0}(\lambda)|=8$, $|\phi_{x_1}(\alpha)|=2$ and
$|\phi_{x_1}(\lambda)|=8$. Since $\phi_{x_i}$ are homomorphisms,
we obtain that $|\phi_{x_{i}}(\gamma)|=|\pm 2p \pm 8q|$, $i=1,2$.
Therefore, we obtain that either
\[
  |\phi_{x_{0}}(\gamma)|+|\phi_{x_{1}}(\gamma)|=|2p+8q|+|2p+8q|,
\]
or
\[
  |\phi_{x_{0}}(\gamma)|+|\phi_{x_{1}}(\gamma)|=|2p+8q|+|2p-8q|.
\]
We claim that the first case can not happen. Suppose not. Take
$\gamma=-4\alpha+\lambda$, then
$|\gamma|=|\phi_{x_{0}}(\gamma)|+|\phi_{x_{1}}(\gamma)|=0$. Hence
$\text{deg}f_{\gamma}=|\gamma|=0$. This is impossible because
$f_{\gamma}$ is not constant. Thus the claim is proved and
$|\gamma|=|(p,q)|=2(|p+4q|+|p-4q|)$.
\end{proof}

Next let $M(0)$ be the closed $3$-manifold obtained by the Dehn
surgery of $M$ along the longitude $\lambda$. It is known that
$M(0)$ admits an essential torus \cite[Page~200]{Boy}. By Theorem
\ref{8res}, we show that $M(0)$ has non-abelian infinite
fundamental group.
\begin{proposition}\label{0sur}
The fundamental group $\pi_1(M(0))$ is a non-abelian, infinite
group.
\end{proposition}

\begin{proof}
Suppose that $\rho\in R(M)$ with the property that its trace
$\sigma(\rho(\lambda))=2$. Now by Proposition \ref{8res}, we have:
\[
  x^4-5x^2+2=2
\]
where $x=\sigma(\rho(\alpha))$. Solve this equation, we get $x=0$
or $x=\pm \sqrt{5}$. By Proposition \ref{8cha}, there exists
$\rho_0\in R(M)$, such that $\sigma(\rho_0(\alpha))=0$ and
$\sigma(\rho_0(\lambda))=2$. In particular the eigenvalues of
$\rho_0(\alpha)$ is $\pm i$. By \cite[Proposition~6.2]{CCGLS},
$\rho_0$ is an irreducible representation. On the other hand,
since $\rho_0(\alpha)$ and $\rho_0(\lambda)$ commute and
$\rho_0(\alpha)$ is not parabolic, $\rho_0(\lambda)$ must be the
identity matrix. Hence $\rho_0$ induce a representation of
$\pi_1(M(0))$. The irreducibility of $\rho_0$ implies that
$\pi_1(M(0))$ is not abelian.

Similarly, let $\rho_{\sqrt{5}}\in R(M)$, such that
$\sigma(\rho_{\sqrt{5}}(\alpha))=\sqrt{5}$ and
$\sigma(\rho_{\sqrt{5}}(\lambda))=2$. By the same reason,
$\rho_{\sqrt{5}}(\lambda)$ equals the identity matrix. Hence
$\rho_{\sqrt{5}}$ induces a representation of $\pi_1(M(0))$. We
can check that the image of $\rho_{\sqrt{5}}$ in
$SL_2(\mathbb{C})$ is torsion-free. Therefore, $\pi_1(M(0))$ is
not finite.
\end{proof}

\begin{remark} We know that $M(0)$ is not a hyperbolic
manifold. $0/1$ is the one of the ten exceptional surgery slopes
of the figure-eight knot. It is interesting to know that we can
prove its fundamental group is non-cyclic, non-abelian and
infinite just from the elementary computations.
\end{remark}

Let $M(3)$ be the closed $3$-manifold obtained by the Dehn surgery
of $M$ along the simple closed curve $\gamma=3\alpha+\lambda$.
\begin{lemma}\label{3/1-irre}
  $M(3)$ has exactly three irreducible $SL_2(\mathbb{C})$ characters.
\end{lemma}
\begin{proof}
 Since $\pi_1(M(3))=\langle \pi_1(M)|\alpha^3\lambda=1\rangle$,
 we have an embedding of character varieties
 $X(M(3))\hookrightarrow X(M)$. On the other hand, $X_0$ is the
 only component of $X(M)$ containing characters of non-abelian
 representations. Thus, the irreducible characters of $M(3)$ are
 contained in the set $S=\{\chi|\chi\in X_0,
 \chi(\alpha^3\lambda)=2\}$. By (\ref{tf3}), we have
 \begin{equation}
  \chi(\alpha^2\lambda)=\chi(\alpha)\chi(\alpha\lambda)-\chi(\lambda),
 \end{equation}
and
  \begin{equation}\label{3/1-t}
    \begin{aligned}
      \chi(\alpha^3\lambda)=&\chi(\alpha)\chi(\alpha^2\lambda)-\chi(\alpha\lambda)\\
                           =&(\chi(\alpha)^2-1)\chi(\alpha\lambda)-\chi(\alpha)\chi(\lambda).
    \end{aligned}
  \end{equation}
By Theorem \ref{8res}, we obtain
 \begin{equation}\label{3/1-eq1}
  (x^2-1)[(4x-x^3)z+x^5-4x^3-x]-x(x^4-5x^2+2)=2.
 \end{equation}
It is clear that the set $S$ is exactly the common solutions to
equations (\ref{3/1-eq1}) and (\ref{8d3}). Note when $x=1$,
(\ref{3/1-eq1}) holds and is independent of the values of $z$.
From (\ref{8d3}), when $x=1$, $z=1$. So $(1,1)\in S$.

We solve $z$ in terms of $x$ from (\ref{3/1-eq1}), then plug in
(\ref{8d3}) and simplify the expression, we get
\begin{equation}\label{3/1-eq2}
  C(x)=x^4-4x^3+2x^2+4x+1=(x^2-2x-1)^2=0.
\end{equation}
It has two solutions $x=1\pm \sqrt{2}$. Hence the set $S$ has
three elements. They are not equal to $\pm 2$ or $\pm \sqrt{5}$.
By Corollary \ref{8-cor2}, $x\ne \pm \sqrt{5}$ implies that each
one is an irreducible character. By Proposition \ref{8cha-1} and
\ref{kin}, $x\ne \pm 2$ means that each one is also a character of
$M(3)$. The result follows.
\end{proof}

Let $sl_2(\mathbb{C})$ be the Lie algebra of $SL_2(\mathbb{C})$.
Then we have the adjoint representation $Ad:
SL_2(\mathbb{C})\rightarrow Aut(sl_2(\mathbb{C}))$. For a
representation $\rho:\pi_1(M(3))\rightarrow SL_2(\mathbb{C})$, let
$H^1(M(3); sl_2(\mathbb{C})_\rho)$ be the first cohomology group
with coefficients in $sl_2(\mathbb{C})$ twisted by the composition
$Ad\circ \rho$.

\begin{proposition}\label{3/1-prop1}
Suppose that $\rho:\pi_1(M(3))\rightarrow SL_2(\mathbb{C})$ is
irreducible. \\ Then $H^1(M(3);sl_2(\mathbb{C})_\rho)=0$.
\end{proposition}

\begin{proof}
By \cite[Theorem~1.1]{BW}, $M(3)$ is not toroidal. Since
$H_1(M(3);\mathbb{Z})$ is finite, $M(3)$ is a small Seifeit
fibered space. By the preceding Lemma \ref{3/1-irre}, $M(3)$ has
irreducible representations, so $\pi_1(M(3))$ is not cyclic. By
\cite[Proposition~7]{BZ1}, $H^1(M(3); sl_2(\mathbb{C})_\rho)=0$.
\end{proof}

\begin{remark} For the other eight exceptional surgery slopes $\pm 1$, $\pm 2$,
$\pm 3$, $\pm 4$, we also have the explicit irreducible
representations with infinite images. Hence their fundamental
groups are all non-abelian and infinite. We omit the details here.
Hence, Proposition \ref{3/1-prop1} holds also for $M(\pm1)$,
$M(\pm2)$ and $M(-3)$.
\end{remark}

Now we turn to the proof of Theorem \ref{8res}. We need to find
explicit expressions of the traces $F(x,z)$, $G(x,z)$ of
$\rho(\lambda)$ and $\rho(\alpha\lambda)$ respectively in terms of
traces of $\rho(\alpha)$ and $\rho(\alpha\beta)$.

We need some preliminary lemmas. For $A$, $B$, $C\in
SL_2(\mathbb{C}))$, we have the following identities on their
traces \cite{Wh}:
\begin{equation}\label{tf1}
  \sigma(AB)=\sigma(BA);
\end{equation}
\begin{equation}\label{tf2}
 \sigma(A)=\sigma(A^{-1});
\end{equation}
\begin{equation}\label{tf3}
  \sigma(AB)=\sigma(A)\sigma(B)-\sigma(AB^{-1});
\end{equation}
\begin{equation}\label{tf4}
 \sigma(ABC)=\sigma(A)\sigma(BC)+\sigma(B)\sigma(AC)+\sigma(C)\sigma(AB)-\sigma(A)\sigma(B)\sigma(C)-\sigma(ACB);
\end{equation}
\begin{equation}\label{tf5}
 \text{For}\, m\geq 2,
 \sigma(A^m)=\sigma(A^{m-1})\sigma(A)-\sigma(A^{m-2}).
\end{equation}

Notice that (\ref{tf1}) is true for any $n\times n$ matrices. It
is easy to see that (\ref{tf4}) and (\ref{tf5}) follow from
(\ref{tf3}). Equations (\ref{tf2}) and (\ref{tf3}) can be checked
by direct
computations.\\

Now for $\rho\in R(M)$, set $a=\rho(\alpha)$, $b=\rho(\beta)$,
$x=\sigma(a)=\sigma(b)$ and
$z=\sigma(ab)=\sigma(\rho(\alpha\beta))$.

\begin{lemma}\label{tr1}
(i) $\sigma(a^2)=x^2-2$.\\
(ii)
$\sigma(a^{-1}b)=\sigma(ba^{-1})=\sigma(b^{-1}a)=\sigma(ab^{-1})=x^2-z$.\\
(iii)
$\sigma(b^{-1}aba^{-1})=\sigma(a^{-1}bab^{-1})=z^2-x^2z+2x^2-2$.
\end{lemma}

\begin{proof}
(i) By (\ref{tf5}), $\sigma(a^2)=\sigma(a)^2-\sigma(I)$, where $I$
is the $2\times 2$ identity matrix.

(ii) By (\ref{tf1}), $\sigma(a^{-1}b)=\sigma(ba^{-1})$ and
$\sigma(b^{-1}a)=\sigma(ab^{-1})$ ; by (\ref{tf2}),
$\sigma(ab^{-1})=\sigma(ba^{-1})$; by (\ref{tf3}),
$\sigma(ab^{-1})=x^2-z$.

(iii) By the subdivision $(b^{-1}a)ba^{-1}$ and (\ref{tf4}), we
have
\[
  \sigma(b^{-1}aba^{-1})=\sigma(b^{-1}a)\sigma(ba^{-1})+\sigma(b)\sigma(b^{-1})
  +\sigma(a^{-1})\sigma(b^{-1}ab)-\sigma(b^{-1}a)\sigma(b)\sigma(a^{-1})-\sigma(I).
\]
Therefore,
\[
  \sigma(b^{-1}aba^{-1})=(x^2-z)^2+x^2+x^2-(x^2-z)x^2-2=z^2-x^2z+2x^2-2.
\]
The proof for $\sigma(a^{-1}bab^{-1})$ is the same by the
subdivision $(a^{-1}b)ab^{-1}$, and we omit it.
\end{proof}

\begin{lemma}\label{tr2}
(i) $\sigma(a^{-1}a^{-1}b)=x(x^2-z)-x$;\\
(ii) $\sigma(ba^{-1}a^{-1}b)=x^4-zx^2-2x^2+2$;\\
(iii) $\sigma(ba^{-1}a^{-1}bab^{-1})=x^2-z$.
\end{lemma}

\begin{proof}
(i) By (\ref{tf3}), we have
\[
  \sigma(a^{-1}(a^{-1}b))=\sigma(a^{-1})\sigma(a^{-1}b)-\sigma(a^{-1}b^{-1}a)
\]
By Lemma \ref{tr1}, each term of the right-hand side is known.
Hence,
\[
  \sigma(a^{-1}a^{-1}b)=x(x^2-z)-x.
\]

(ii) By (\ref{tf3}),
\[
  \sigma(b(a^{-1}a^{-1}b))=\sigma(b)\sigma(a^{-1}a^{-1}b)-\sigma(bb^{-1}aa)=\sigma(b)\sigma(a^{-1}a^{-1}b)-\sigma(a^2)
\]
The formula follows.

(iii) By (\ref{tf4}), we have
\begin{align*}
  \begin{split}
    \sigma((ba^{-1})(a^{-1}b)(ab^{-1}))&=\sigma(ba^{-1})\sigma(a^{-1}bab^{-1})+\sigma(a^{-1}b)\sigma(1)+\sigma(ab^{-1})\sigma(ba^{-1}a^{-1}b)-\\
                                       & \sigma(ba^{-1})\sigma(a^{-1}b)\sigma(ab^{-1})-\sigma(a^{-1}b)
  \end{split}
\end{align*}

Plug in what we know on the right-hand side and simplify, we
obtain the formula.
\end{proof}

\begin{lemma}\label{tr3}
(i) $\sigma(ab^{-1}a)=\sigma(aab^{-1})=x^3-zx-x$;\\
(ii) $\sigma(aab)=xz-x$;\\
(iii) $\sigma(ab^{-1}aba^{-1})=x$
\end{lemma}

\begin{proof}
(i) $\sigma(a(b^{-1}a))=\sigma(a)\sigma(b^{-1}a)-\sigma(b)$, and $\sigma(a(ab^{-1}))=\sigma(a)\sigma(ab^{-1})-\sigma(aba^{-1})$;\\
(ii) $\sigma(a(ab))=\sigma(a)\sigma(ab)-\sigma(ab^{-1}a^{-1})$\\
(iii)
$\sigma(a(b^{-1}aba^{-1}))=\sigma((b^{-1}aba^{-1})a)=\sigma(b^{-1}ab)=x$.
\end{proof}

\begin{proof}[Proof of Theorem~\ref{8res}]
First, let us compute $F(x,z)=\sigma(\rho(\lambda))$, the trace of
$\rho(\lambda)$. By (\ref{tf4}), we have
\begin{align*}
  \begin{split}
     \sigma((b^{-1}a)(ba^{-1})(a^{-1}bab^{-1}))&=\sigma(b^{-1}a)\sigma(ba^{-1}a^{-1}bab^{-1})+\sigma(ba^{-1})\sigma(ab^{-1})+\\
                                               & \sigma(a^{-1}bab^{-1})\sigma(b^{-1}aba^{-1})
                   -\sigma(b^{-1}a)\sigma(ba^{-1})\sigma(a^{-1}bab^{-1})-\sigma(I)
  \end{split}
\end{align*}

By Lemmas \ref{tr1} and \ref{tr2}, we know all the terms on the
right-hand side. Notice that $(x,z)\in X_0$, hence
$z^2-(1+x^2)z+2x^2-1=0$. Now divide the right-hand side by
$z^2-(1+x^2)z+2x^2-1$, the remainder is $F(x,z)$. we calculate that $F(x,z)=x^4-5x^2+2$.\\

For the trace of $\rho(\alpha\lambda)$, we have
\begin{align*}
  \begin{split}
     \sigma((ab^{-1}a)(ba^{-1})(a^{-1}bab^{-1}))&=\sigma(ab^{-1}a)\sigma(ba^{-1}a^{-1}bab^{-1})+\sigma(ba^{-1})\sigma(aab^{-1})+\\
                                                & \sigma(a^{-1}bab^{-1})\sigma(ab^{-1}aba^{-1})-\sigma(ab^{-1}a)\sigma(ba^{-1})\sigma(a^{-1}bab^{-1})
                                                -\sigma(a).
  \end{split}
\end{align*}

By Lemmas  \ref{tr1}, \ref{tr2} and \ref{tr3}, we have all the
terms on the right-hand side. Then we divide the result of the
right-hand side by the polynomial $z^2-(1+x^2)z+2x^2-1=0$ and
$G(x,z)$ equals the remainder. We calculate that it is
$(4x-x^3)z+(x^5-4x^3-x)$.
\end{proof}

\end{document}